\documentclass[a4paper,12pt]{amsart}
\usepackage{amssymb}
\usepackage{ifthen}
\usepackage{mathrsfs}
\nonstopmode
\numberwithin{equation}{section}
\theoremstyle{plain}
\newtheorem{thm}{Theorem}[section]

\newtheorem{lem}[thm]{Lemma}

\theoremstyle{definition}

\newenvironment{rem}{%
\bigskip
\noindent
\textsl{{\sl Remark. }}}{\bigskip}

\newenvironment{pf}[1][]{%
 \vskip 3mm
 \noindent
 \ifthenelse{\equal{#1}{}}%
  {{\slshape Proof. }}%
  {{\slshape #1.} }%
 }%
{\qed\bigskip}

\newcounter{alphabet}
\newcounter{tmp}
\newenvironment{Thm}[1][]{\refstepcounter{alphabet}%
\bigskip%
\noindent%
{\bf Theorem \Alph{alphabet}}%
\ifthenelse{\equal{#1}{}}{}{ (#1)}%
{\bf .}
\itshape}{\vskip 8pt}

\newcommand{\A}{{\mathscr A}}

\newcommand{\es}{{\mathscr S}}

\newcommand{\R}{{\mathbb R}}
\newcommand{\C}{{\mathbb C}}

\newcommand{\D}{{\mathbb D}}
\newcommand{\sphere}{{\widehat{\mathbb C}}}

\renewcommand{\Re}{{\operatorname{Re}\,}}

\renewcommand{\mod}{{\operatorname{mod}\,}}

\renewcommand{\arg}{\,{\operatorname{arg}\,}}
\newcommand{\argl}{\,{\operatorname{arg}_\lambda\,}}
\newcommand{\Sp}{{\mathscr{F}}}

\newcounter{minutes}
\divide\time by 60
\newcounter{hours}
\multiply\time by 60
\addtocounter{minutes}{-\time}

\begin{document}

\title[Spirallike functions and starlike functions]{
Correspondence between spirallike functions and starlike functions }

\author[Y.~C.~Kim]{Yong Chan Kim}
\address{Department of Mathematics Education, Yeungnam University, 214-1 Daedong
Gyongsan 712-749, Korea}
\email{kimyc@ynu.ac.kr}
\author[T. Sugawa]{Toshiyuki Sugawa}
\address{Graduate School of Information Sciences,
Tohoku University, Aoba-ku, Sendai 980-8579, Japan}
\email{sugawa@math.is.tohoku.ac.jp}
\begin{abstract}
Let $\lambda$ be a real number with $-\pi/2<\lambda<\pi/2.$ In order
to study $\lambda$-spirallike functions, it is natural to measure
the angle according to $\lambda$-spirals. Thus we are led to the
notion of $\lambda$-argument. This fits well the classical
correspondence between $\lambda$-spirallike functions and starlike
functions. Using this idea, we extend deep results of Pommerenke and
Sheil-Small on starlike functions to spirallike functions. As an
application, we solved a problem given by Hansen in \cite{Hansen71}.
\end{abstract}
\thanks{
The first author was supported by Yeungnam University.
The second author was supported in part by JSPS Grant-in-Aid for
Exploratory Research, 19654027.
}
\maketitle

\section{Introduction}
A domain $\Omega$ with $0\in\Omega$ in the complex plane $\C$ is called
{\it starlike} with respect to $0$
if the line segment $[0,w]$ joining $0$ and $w$
is contained in $\Omega$ whenever $w\in\Omega.$
Note that a starlike domain is simply connected.
An analytic function $f$ on the unit disk $\D=\{z\in\C: |z|<1\}$
with $f(0)=0$
is called {\it starlike} if $f$ maps $\D$ univalently onto a starlike domain
with respect to $0.$
It is well known that starlikeness is characterized by the condition
$$\Re(zf'(z)/f(z))>0,~|z|<1.$$
Starlike functions have been studied by many authors.
See, for instance, Duren \cite{Duren:univ} and Goodman \cite{Goodman:univ}
and references therein.

The notion of starlike domains and starlike functions can be extended
by using logarithmic spirals instead of line segments.
Let $\lambda$ be a real number between $-\pi/2$ and $\pi/2.$
The curve
$$\gamma_\lambda:t\mapsto\exp(te^{i\lambda}),~t\in\R,$$
and their rotations $e^{i\theta}\gamma_\lambda,~\theta\in\R,$
are called $\lambda$-spirals.
These curves $\gamma(t)$
are characterized (up to parametrization)
by the property that the oriented angle from $\gamma(t)$ to the tangent
vector $\gamma'(t),$  which is called the {\it radial angle}, is constantly
$\lambda;$ in other words,
$$\arg(\gamma'(t)/\gamma(t))=\lambda.$$
Also, note that this curve family is invariant under the dilation
$z\mapsto cz$ for $c\in\C\setminus\{0\}.$

For $w\in\C,$ we define the $\lambda$-spiral segment $[0,w]_\lambda$
by
$$
[0,w]_\lambda=w\cdot \gamma_\lambda((-\infty,0])\cup\{0\}
=\{w\exp(te^{i\lambda}): t\le0\}\cup\{0\}.
$$
Clearly, $[0,w]_0$ is the line segment $[0,w].$
A domain $\Omega$ with $0\in\Omega$ is called {\it $\lambda$-spirallike}
(with respect to $0$)
if $[0,w]_\lambda\subset \Omega$ whenever $w\in\Omega.$
Similarly, an analytic function $f$ on the unit disk $\D$
with $f(0)=0$
is called {\it $\lambda$-spirallike} if $f$ maps $\D$ univalently onto
a $\lambda$-spirallike domain.
It is also known that an analytic function $f$ on $\D$ with $f(0)=0$
is $\lambda$-spirallike if and only if
\begin{equation}\label{eq:sp}
\Re\left(e^{-i\lambda}\frac{zf'(z)}{f(z)}\right)>0,\quad |z|<1.
\end{equation}
For the proof and a geometric interpretation of this condition,
the reader should consult \S 2.7 of Duren's book \cite{Duren:univ}.
We remark that many authors adopt the condition
with $i\lambda$ instead of $-i\lambda$ in \eqref{eq:sp} as
the definition of $\lambda$-spirallike functions
(for inctance, see \cite[\S 9.3]{Goodman:univ}).

Without much loss of generality, we may focus on analytic functions
$f$ on $\D$ with $f(0)=0, f'(0)=1$ in the sequel, and we denote by
$\A$ the set of such normalized functions. Let $\es$ stand for the
set of univalent functions in $\A.$ We further denote by
$\Sp_\lambda$ the subclass of $\A$ consisting of functions
satisfying \eqref{eq:sp} for $-\pi/2<\lambda<\pi/2.$ In particular,
$\Sp_0$ coincides with the class $\es^*$ of starlike functions in
$\A.$

A function in the union $\Sp=\bigcup_{|\lambda|<\pi/2}\Sp_\lambda$
is simply called {\it spirallike}. Note that $\Sp\subset\es.$ It is
not necessarily easy to deal with spirallike functions in spite of
its simple form of the characterizing condition in \eqref{eq:sp}.
For instance, $\Sp_\lambda$ is not contained in the class of
close-to-convex functions for $\lambda\ne0$ (see \cite[\S
2.7]{Duren:univ}). Therefore, a relatively small number of papers
have been devoted to the study of spirallike functions so far.

For starlike functions $f,$ it is fundamental to consider the radial limit
of the argument (see \cite{Pom62} for instance):
$$\arg f(e^{it})=\displaystyle\lim_{r\to1-}\arg f(re^{it}).$$
However, it is not appropriate to consider the same quantity
for spirallike functions because the limit might not exist.
It is more natural to measure the angle by using logarithmic spirals.

In the present paper, we propose the notion of {\it
$\lambda$-argument} (cf.~\cite{Yos71}), denoted by $\argl.$ We set
$$\theta=\argl w \quad\text{if}\quad w\in e^{i\theta}\gamma_\lambda(\R).$$
Note that arg$_0 w=\arg w.$ Note also that we have a freedom for the
choice of $\argl w$ up to an integer multiple of $2\pi$ as in the
case of $\arg w.$ Elimination of this ambiguity for certain cases
and further properties of $\lambda$-argument will be discussed in
Section 2. By means of $\lambda$-argument, we can state one of our
main results as the following form, which is a straightforward
generalization of \cite[Theorem 3.18]{Pom:bound} on starlike
functions (originally due to Pommerenke \cite{Pom62}, \cite{Pom63}
and Sheil-Small \cite{Sheil70}).

\begin{thm}[Representation Theorem]\label{thm:main1}
Let $f\in\Sp_\lambda$ for a $\lambda\in(-\pi/2,\pi/2).$
Then the limits
$$
\beta(t)=\lim_{r\to1-}\argl f(re^{it})
\quad\text{and}\quad
f(e^{it})=\lim_{r\to1-} f(re^{it})\in\sphere=\C\cup\{\infty\}
$$
exist for every $t\in\R$ in such a way that
$\beta(t)$ is non-decreasing in $t$ and $\beta(t+2\pi)=\beta(t)+2\pi.$
Moreover, $f$ is represented by
\begin{equation}\label{eq:f}
f(z)=z\exp\left(
-\frac{e^{i\lambda}\cos\lambda}\pi\int_0^{2\pi}\log(1-e^{-it}z)d\beta(t)
\right),\quad z\in\D.
\end{equation}
Conversely, if $\beta(t)$ is a non-decreasing real-valued function in $t\in\R$
with $\beta(t+2\pi)=\beta(t)+2\pi,$ then the function $f$ given by \eqref{eq:f}
is $\lambda$-spirallike.
\end{thm}

We remark that the representation formula \eqref{eq:f} itself is not new
(cf.~formula (14) in MacGregor \cite{Mac73}) and indeed it immediately follows
from the Herglotz formula.
We emphasize that we now have a geometric interpretation of the measure
$d\beta(t)$ in the representation formula.

Let $f\in\es.$
For $R>0,$ we denote by $\alpha(R,f)$ the length of the largest arc
contained in the set $\{\zeta\in\partial\D: R\zeta\in f(\D)\}.$
By definition, we have $0\le \alpha(R,f)\le 2\pi$ for $R>0.$
Obviously, $\alpha(R,f)$ is non-increasing in $R$ for a spirallike function $f.$
Therefore, the limit
$$
A(f)=\lim_{R\to+\infty}\alpha(R,f)
$$
exists and satisfies $0\le A(f)\le 2\pi$ for $f\in\Sp.$
Let
$$
M(r,f)=\max_{|z|=r}|f(z)|.
$$
and let $q_0=\frac1\pi A(f)\cos^2\lambda$ for $f\in\Sp_\lambda.$
Hansen \cite{Hansen71} showed that
$$\sup\{p>0: f\in H^p\}=1/q_0$$
and, as an application of this result, that
$$M(r,f)=o[(1-r)^{-q}],~q>q_0,$$
for $f\in\Sp_\lambda.$ In other words,
$$\limsup_{r\to1-}\log M(r,f)/\log\frac1{1-r}\le q_0.$$
We will show
the following refinement. Note that Pommerenke \cite{Pom62} proved
it for starlike functions (see the proof of Theorem \ref{thm:main2}
for details).

\begin{thm}\label{thm:main2}
Let $f\in\Sp_\lambda.$
Then
$$
\lim_{r\to1-} \frac{\log M(r,f)}{\log\frac1{1-r}}=\frac{A(f)\cos^2\lambda}\pi.
$$
\end{thm}

He suspected in \cite{Hansen71}
that

$$M(r,f)=O[(1-r)^{-q_0}]$$

if $A(f)\ne0.$ We will show that this is not true in general.

\begin{thm}\label{thm:ce}
Let $\lambda\in(-\pi/2,\pi/2)$ and $0<A<2\pi.$
Then there is a $\lambda$-spirallike function $f$ with
$A(f)=A$ so that $M(r,f)=O[(1-r)^{-A(f)\cos^2\lambda/\pi}]$
does not hold.
\end{thm}

\section{Preliminaries}

We first summarize basic properties of the $\lambda$-argument.
The following elementary lemma is convenient in various computations.

\begin{lem}\label{lem:ele}
For $\lambda\in(-\pi/2,\pi/2),~\theta\in\R$ and $w\in\C\setminus\{0\},$
$$
\argl w=\arg w-(\tan\lambda)\log|w|  \ \ (\mod 2\pi).
$$
\end{lem}

\begin{pf}
Let $\theta=\argl w.$
Then, by definition,
$w=e^{i\theta}\gamma_\lambda(t)=\exp(i\theta+te^{i\lambda})$
for some $t\in\R.$
This is equivalent to the relation $\log w=i\theta+te^{i\lambda}\
(\mod 2\pi i).$
Taking the real and imaginary part, we have
$\log|w|=t\cos\lambda$ and $\arg w=\theta+t\sin\lambda~(\mod 2\pi).$
We now eliminate $t$ from these two relations to obtain the required one.
We can trace back the above procedure to check the converse.
\end{pf}

With the help of the last lemma, we can easily check the following
analog to a familiar law for the ordinary argument.

\begin{lem}
For nonzero complex numbers $w_1, w_2$ and $\lambda\in(-\pi/2,\pi/2),$
$$
\argl(w_1w_2)=\argl w_1+\argl w_2
\quad (\mod 2\pi).
$$
\end{lem}

Also, by Lemma \ref{lem:ele}, we have a canonical way to take a harmonic
branch of $\argl h$ for a non-vanishing analytic function $h.$

\begin{lem}
Let $\lambda\in(-\pi/2,\pi/2)$ and
let $h$ be a non-vanishing analytic function in the unit disk $\D$
with $h(0)=1.$
Then there is a unique harmonic function $u$ on $\D$ with $u(0)=0$ such that
$\argl h(z)=u(z) ~(\mod 2\pi)$ for $z\in\D.$
\end{lem}

\begin{pf}
Since $h$ is non-vanishing on the simply connected domain $\D,$
we can take a harmonic branch $v$ of $\arg h$ on $\D$ so that $v(0)=0.$
Then the harmonic function $u=v-(\tan\lambda)\log|h|$ satisfies
the required conditions.
\end{pf}

In what follows, unless otherwise stated, we will take the above $u$ as a branch
of $\argl h$ for a non-vanishing analytic function $h$ on $\D$
with $h(0)=1.$
The same applies to the ordinary argument.
For instance, we can take $f(z)/z$ as $h(z)$ for $f\in\es.$
When $f$ is starlike, the following deep properties were proved by
Pommerenke \cite{Pom62} and Sheil-Small \cite{Sheil70}
(see also \cite[\S 3.6]{Pom:bound}).


\begin{Thm}\label{Thm:PS}
Let $g\in\es^*.$
Then the limits
$$
U(t)=\lim_{r\to1-}\arg\frac{g(re^{it})}{re^{it}}
\quad\text{and}\quad
g(e^{it})=\lim_{r\to1-} g(re^{it})\in\sphere
$$
exist for every $t\in\R,$ and
$\beta(t,g)=U(t)+t$ is a non-decreasing function in $t$ with
$\beta(t+2\pi,g)=\beta(t,g)+2\pi.$
Moreover, the left and right limits of $\beta(t,g)$ satisfy
the following relation:
\begin{equation}\label{eq:rl}
\beta(t,g)=\frac12\big(\beta(t+,g)+\beta(t-,g)\big).
\end{equation}
For $t_0\in\R$ and $\theta\in(0,2\pi],$ the relation
$\beta(t_0+,g)-\beta(t_0-,g)=\theta$ holds if and only if
the image domain $g(\D)$ contains a maximal sector of the form
$\{w: |\arg w-\beta(t_0,g)|<\theta/2\}.$
\end{Thm}

We denote by $B(g)$ the maximal jump of $\beta(t,g)$ for $g\in\es^*.$
In other words,
$$
B(g)=\max_{t\in\R}\big[\beta(t+,g)-\beta(t-,g)\big].
$$
Then, as an immediate consequence of Theorem A, 
we have the relation
\begin{equation}\label{eq:AB}
A(g)=B(g) \qquad \text{for}\quad g\in\es^*.
\end{equation}

As for the quantity $B(g),$ Pommerenke \cite{Pom62} found
a connection with the growth of a starlike function $g.$

\begin{Thm}[Pommerenke \cite{Pom62}]\label{Thm:growth}
Let $g\in\es^*.$
Then
$$
\lim_{r\to1-} \frac{\log M(r,g)}{\log\frac1{1-r}}=\frac{B(g)}\pi.
$$
\end{Thm}

\section{Spirallike counterpart}

In this section, we extend the results of Pommerenke and Sheil-Small
in the previous section to spirallike functions.
Our main tools will be the $\lambda$-argument and a
useful correspondence between $\lambda$-spirallike functions
and starlike functions.

Let $f\in\Sp_\lambda.$ Then, by \eqref{eq:sp}, we find an analytic
function $p$ on $\D$ with $\Re p>0$ and $p(0)=1$ such that
$$
e^{-i\lambda}\frac{zf'(z)}{f(z)}
=p(z)\cos\lambda-i\sin\lambda,
\quad z\in\D.
$$
Define $g\in\es^*$ by the relation $zg'(z)/g(z)=p(z).$
The correspondence $f\mapsto g$ gives a bijection from
$\Sp_\lambda$ onto $\es^*=\Sp_0.$
Integrating the above relation,
we arrive at the following well-known fact (cf.\ \cite{BK70}).

\begin{lem}\label{lem:fg}
Let $\lambda\in(-\pi/2,\pi/2).$
There corresponds to $f\in\Sp_\lambda$ a unique starlike function
$g\in\es^*$ in such a way that
\begin{equation}\label{eq:fg}
\frac{f(z)}{z}=\left(\frac{g(z)}{z}\right)^{e^{i\lambda}\cos\lambda},
\quad z\in\D.
\end{equation}
\end{lem}

This relation serves as a key to reduce a problem concerned with
spirallike functions to one with starlike functions.
Note here, however, that the relation \eqref{eq:fg} does not give
a transformation of the image domain $f(\D)$ onto $g(\D)$
because the term $z$ is involved.

To realize the connection, we extend notions for starlike functions
to spirallike ones by using the $\lambda$-argument. Let
$\lambda\in(-\pi/2,\pi/2).$ We call the set
$$
S_\lambda(\theta_0,\alpha)=
\{ e^{i\theta}\gamma_\lambda(t):
t\in\R, |\theta-\theta_0|<\alpha/2\}
$$
a $\lambda$-spiral sector of opening $\alpha$ with
center angle $\theta_0.$
Here, recall that $\gamma_\lambda(t)=\exp(e^{i\lambda}t).$
A $\lambda$-spiral sector $S$ of opening $\alpha$ is said to be
maximal in a domain $\Omega$ if $S\subset\Omega$ and
if there are no $\alpha'>\alpha$ and $\theta_0\in\R$
such that $S\subset S_\lambda(\theta_0,\alpha')\subset\Omega.$

We first extend Theorem A 
to spirallike functions.

\begin{thm}\label{thm:spiral}
Let $f\in\Sp_\lambda$ for a $\lambda$ with $-\pi/2<\lambda<\pi/2.$
Then the limits
$$
U_\lambda(t)=\lim_{r\to1-}\argl\frac{f(re^{it})}{re^{it}}
\quad\text{and}\quad
f(e^{it})=\lim_{r\to1-} f(re^{it})\in\sphere
$$
exist for every $t\in\R,$ and
$\beta_\lambda(t,f)=U_\lambda(t)+t$ is a non-decreasing function in $t$ with
$\beta_\lambda(t+2\pi,f)=\beta_\lambda(t,f)+2\pi.$
Moreover, the left and right limits of $\beta_\lambda(t,f)$ satisfy
the following relation:
\begin{equation}\label{eq:rl2}
\beta_\lambda(t,f)=\frac12\big(\beta_\lambda(t+,f)+\beta_\lambda(t-,f)\big).
\end{equation}
For $t_0\in\R$ and $\theta\in(0,2\pi],$ the relation
$\beta_\lambda(t_0+,f)-\beta_\lambda(t_0-,f)=\theta$ holds if and only if
the image domain $f(\D)$ contains a maximal $\lambda$-spiral sector
of the form $S_\lambda(\beta_\lambda(t_0,f),\theta).$
\end{thm}

\begin{pf}
For $f\in\Sp_\lambda,$ we take $g\in\es^*$ as in Lemma \ref{lem:fg}.
Then, by taking the imaginary part of the relation
$$
\log\frac{g(z)}{z}=\frac{e^{-i\lambda}}{\cos\lambda}\log\frac{f(z)}{z},
$$
we obtain
$$
\arg\frac{g(z)}{z}=\arg\frac{f(z)}{z}
-(\tan\lambda)\log\left|\frac{f(z)}{z}\right|.
$$
By Lemma \ref{lem:ele}, we now have the useful formula
\begin{equation}
\arg\frac{g(z)}{z}=\argl\frac{f(z)}{z}.
\end{equation}
It now follows from Theorem A 
that the limits $U_\lambda(t)$ and $f(e^{it})$ exist
and $\beta_\lambda(t,f)$ satisfies the relation
\begin{equation}\label{eq:beta}
\beta_\lambda(t,f)=\beta(t,g),\quad t\in\R.
\end{equation}
In this way, all the assertions in the theorem but the last can be checked.
To show the last assertion, we assume that
$$\beta_\lambda(t_0+,f)-\beta_\lambda(t_0-,f)=\theta>0$$
and let
$$\theta_0=\beta_\lambda(t_0,f).$$
We first check that $S=S_\lambda(\theta_0,\theta)$
is contained in $f(\D).$
Suppose, to the contrary, that there is a point $w_0\in S\setminus f(\D).$
Let $\theta_1=\argl w_0$ and let $w_1$ be the other end point of the
curve $e^{i\theta_1}\gamma_\lambda(\R)\cap f(\D)$ than $0.$
Note that $w_1\ne\infty.$
Since $w_1$ is accessible along the curve $e^{i\theta_1}\gamma_\lambda$
in $f(\D),$ there exists a $t_1\in\R$ such that $f(e^{it_1})=w_1.$
By adding an integer multiple of $2\pi$ to $\theta_1$ if necessary,
we may assume that $\beta_\lambda(t_1,f)=\theta_1.$
Now the inequalities
$$\beta_\lambda(t_0-,f)<\beta_\lambda(t_1,f)
=\theta_1<\beta_\lambda(t_0+,f)$$
enforce $t_1=t_0$ and thus
$\theta_1=\theta_0.$ In this way, we see that $S\setminus f(\D)$ is
contained in the $\lambda$-spiral $e^{i\theta_0}\gamma_\lambda(\R).$
Since $f(\D)$ is simply connected, $S\setminus f(\D)$ must be a
closed $\lambda$-spiral ray with tip at $f(e^{it_0})\ne\infty.$ In
particular, $f(e^{it})$ must be continuous at $t=t_0,$ which
contradicts the assumption $\theta>0.$ Thus, we have shown that
$$S\subset f(\D).$$
Since there are sequences $t_n^-<t_0<t_n^+$ such
that $t_n^-\to t_0, t_n^+\to t_0$ and that $w_n^-=f(e^{it_n^-})$ and
$w_n^+=f(e^{it_n^+})$ are finite points. Since
$$\argl w_n^\pm\to \beta_\lambda(t_0^\pm,f),$$
one can easily see that $S$ is a maximal $\lambda$-spiral sector
in $f(\D).$

The converse can be checked by using the relation in \eqref{eq:rl2}.
\end{pf}

We now prove Theorem \ref{thm:main1}.

\begin{pf}[Proof of Theorem \ref{thm:main1}]
Note that $\beta(t)$ in the theorem is nothing but $\beta_\lambda(t,f).$
The formula
$$
g(z)=z\exp\left(
-\frac{1}\pi\int_0^{2\pi}\log(1-e^{-it}z)d\beta(t,g)
\right),\quad z\in\D
$$
is well known for a starlike function $g$ (cf.~\cite[Theorem 3.18]{Pom:bound}).
We now use the relations \eqref{eq:fg} and \eqref{eq:beta} to
deduce \eqref{eq:f}.
The other assertions follow from Theorem \ref{thm:spiral}.
\end{pf}

Let
$$
B_\lambda(f)=\max_{t\in\R}\big[\beta_\lambda(t+,f)-\beta_\lambda(t-,f)\big]
$$
for $f\in\Sp_\lambda.$
By Theorem \ref{thm:spiral}, we have the relation
$A(f)=B_\lambda(f).$
On the other hand, by \eqref{eq:beta}, we obtain
$B(g)=B_\lambda(f)$ for $f$ and $g$ in Lemma \ref{lem:fg}.
We now summarize these formulas and \eqref{eq:AB} as in the following.

\begin{lem}\label{lem:same}
Suppose that $f\in\Sp_\lambda$ and $g\in\es^*$ are related by \eqref{eq:fg}.
Then
$$
A(f)=B_\lambda(f)=B(g)=A(g).
$$
\end{lem}

We are ready to prove Theorem \ref{thm:main2}.

\begin{pf}[Proof of Theorem \ref{thm:main2}]
For a given $f\in\Sp_\lambda,$ we take a $g\in\es^*$ as in Lemma \ref{lem:fg}.
Then by \eqref{eq:fg} we have
$$
\log\left|\frac{f(z)}{z}\right|=
\cos^2\lambda\cdot\log\left|\frac{g(z)}{z}\right|
-\sin\lambda\cos\lambda\cdot \arg\left(\frac{g(z)}{z}\right).
$$
Since
$$|\arg(g(z)/z)|\le 2\arcsin |z|<\pi$$
for a starlike function $g$
by a theorem of Goodman \cite[Theorem 2]{Good53},
the second term in the right-hand side is bounded.
This implies
\begin{equation}\label{eq:M}
\log M(r,f)=\cos^2\lambda\cdot \log M(r,g)+O(1) \quad (r\to1-),
\end{equation}
and therefore, by Theorem B, 
$$
\lim_{r\to1-} \frac{\log M(r,f)}{\log\frac1{1-r}}
=\cos^2\lambda\lim_{r\to1-} \frac{\log M(r,g)}{\log\frac1{1-r}}
=\cos^2\lambda\cdot\frac{B(g)}\pi.
$$
In view of Lemma \ref{lem:same}, we now have the required relation.
\end{pf}

\section{Examples}

As is noted in Introduction, Hansen\cite{Hansen71} suspected that
$M(r,f)=O[(1-r)^{-q_0}]$ for $f\in\Sp_\lambda,$ where $q_0=\frac1\pi
A(f)\cos^2\lambda.$ To give the answer negatively, we construct
examples of starlike and spirallike functions. We start with a
simple one. The function $g_0$ below is a standard example
(cf.~\cite[p.~51, Exercise 2]{Duren:hp}). It seems, however, that
attention is not paid to its geometric properties so far.

In what follows, we will mean by
$$\varphi(r)\asymp \psi(r) ~(r\to 1-)$$
that there exist positive constants $A, B$ and $r_0\in(0,1)$ such that
$$A\psi(r)\le \varphi(r)\le B\psi(r)$$
holds for $r_0\le r<1.$

\begin{lem}
Let
$$
g_0(z)=\frac1{1-z}\log\frac1{1-z}
=\sum_{n=1}^\infty\left(1+\frac12+\cdots+\frac1n\right)z^n,
\quad |z|<1.
$$
Then $g_0$ is a starlike function with $A(g_0)=\pi$ and satisfies
$$
M(r,g_0)\asymp \frac{\log\frac1{1-r}}{1-r} \quad (r\to1-).
$$
\end{lem}

\begin{pf}
A simple computation gives
$$
\frac{zg_0'(z)}{g_0(z)}
=\frac{z}{1-z}+G(z),
$$
where
$$
G(z)=-\frac{z}{(1-z)\log(1-z)}.
$$
By a theorem of Wilken and Feng \cite{WF80},
$$\Re G(z)>G(-1)=1/(2\log2)$$
for $z\in\D.$ Therefore,
$$
\Re\frac{zg_0'(z)}{g_0(z)}>-\frac12+\frac1{2\log2}=0.2213\dots>0,\quad
z\in\D.
$$
Thus $g_0$ is starlike.
Since $\arg g_0(e^{i\theta})\to \pm\pi/2$ as $\theta\to 0\pm,$ we can
easily see that $A(g_0)=\pi.$
The asymptotic behavior of $M(r,g_0)$ is obvious by the form of $g_0.$
\end{pf}

Hansen showed also the following in \cite{Hansen71}: Let
$f(z)=z+\sum_{n=2}^\infty a_nz^n$ be in $\Sp_\lambda.$ Then
$a_n\to0~(n\to\infty)$ when $q_0<1,$ and $a_n=o(n^{q-1})$ for any
$q>q_0$ when $q_0\ge1.$ He suspected that
$$a_n=O(n^{q_0-1})$$
if $q_0\ge1.$ Since
$$1+1/2+\dots+1/n\asymp \log n,$$
the above example shows that it does not hold when $\lambda=0$ and
$q_0=1.$

In order to prove Theorem \ref{thm:ce}, we need to modify the above example.
To this end, we prepare some lemmas.

\begin{lem}\label{lem:C}
There exists a number $C_0>2$ such that the following inequalities
hold for $z\in\D$ whenever $C\ge C_0:$
\begin{align*}
&\Re\left[\frac1{(1-z)\log\frac C{1-z}}\right]>\frac1{2\log\frac C2}, \\
&\Re\left[\frac z{(1-z)\log\frac C{1-z}}\right]>\frac{-1}{2\log\frac C2}.
\end{align*}
\end{lem}

\begin{pf}
Set
$$p(z)=1/[(1-z)\log\frac{C}{1-z}]$$
and
$$q(z)=z/[(1-z)\log\frac{C}{1-z}].$$
Assume that $|z|=1$ and $z\ne1.$ Since $\Re[1/(1-z)]=1/2,$ we can
write
$$
\frac1{1-z}=\frac{1+i\tan\theta}2,
$$
where
$$\theta=\arg\frac{1}{1-z}=-\arg(1-z)\in(-\pi/2,\pi/2).$$
By the symmetry, we can assume that $0<\theta<\pi/2.$
An elementary computation yields the expression
$$
\Re p(z)=\frac12\cdot\frac{\log\frac{C}{2\cos\theta}+\theta\tan\theta}%
{(\log\frac{C}{2\cos\theta})^2+\theta^2}.
$$
The assertion
$$\Re p(z)\ge 1/(2\log(C/2))$$
is thus equivalent to the inequality
$$
\log\frac C2\ge \frac{(\log\cos\theta)^2+\theta^2}%
{\theta\tan\theta+\log\cos\theta}=:Q(\theta).
$$
It is easy to see that
$$\lim_{\theta\to0+}Q(\theta)=2 \quad\text{and}\quad
\lim_{\theta\to\pi/2-}Q(\theta)=0.$$
By the continuity of the
function $Q(\theta),$ we now conclude that $Q(\theta)$ is bounded.
Thus we can take a constant $C_0\ge2e^2$ so that
$$Q(\theta)\le \log(C_0/2)$$
holds for $0<\theta<\pi/2.$ Then
$$\Re p(z)\ge 1/(2\log(C/2))$$
holds for $|z|=1, z\ne1$ and $C\ge C_0.$ To conclude that
$$\Re p(z)> 1/(2\log(C/2))$$
for $z\in\D,$ by the minimum principle, it is enough to check the
condition
$$
\liminf_{z\to 1 ~\text{in}~\D}\frac{\Re p(z)}{P_0(z)}\ge0,
$$
where
$$P_0(z)=\frac{1-|z|^2}{|1-z|^2}$$
is the Poisson kernel. This can be confirmed by looking at the
expression
$$
\frac{\Re p(z)}{P_0(z)}=\left(|1-z|\left|\log\tfrac C{1-z}\right|\right)^{-2}
\left[(1-x)\log\tfrac C{|1-z|}-y\arg(1-z)\right]
$$
for $z=x+iy\in\D.$
Here, we use the inequality
$$-y\arg(1-z)\ge0.$$

To prove the second inequality, we express the function $q$ by
$$q(z)=p(z)-1/h(z),$$
where
$$h(z)=\log C/(1-z).$$ Since $1+zh''(z)/h'(z)=1/(1-z)$ has
positive real part, $h$ maps $\D$ onto a convex domain symmetric in
$\R.$ Therefore,
$$h(\D)\subset\{w: \Re w>h(-1)\}.$$
In particular,
$$h(\D)\cap\{w: |w-h(-1)/2|\le h(-1)/2\}=\emptyset,$$
which is equivalent to
$$\Re(1/h(z))<1/h(-1),~z\in\D.$$
Therefore, we have now
$$
\Re q(z)=\Re p(z)-\Re\frac1{h(z)}>p(-1)-\frac1{h(-1)}
=q(-1).
$$
Thus we have shown the second inequality.
\end{pf}

\begin{rem}
It seems that $Q(\theta)$ is monotone decreasing in $0<\theta<\pi/2.$
If this is the case, we have $\sup Q(\theta)=2$ and thus
we can take $2e^2=14.778\dots$ as $C_0$ in the above lemma.
\end{rem}

We are now ready to show the following.

\begin{lem}\label{lem:example}
Let $0<\alpha<2.$
Choose positive numbers $\beta$ and $c$ so that
$$
c\le \frac1{\log C_0}
\quad\text{and}\quad
\alpha+\frac{c\beta}{1-c\log2}<2,
$$
where $C_0$ is the number appearing in Lemma \ref{lem:C}.
Then the function
$$
g(z)=\frac{z}{(1-z)^\alpha}\left(1+c\log\frac1{1-z}\right)^\beta
$$
is starlike and satisfies $A(g)=\pi\alpha$ and
$$M(r,g)\asymp (1-r)^{-\alpha}(\log\frac1{1-r})^\beta$$
as $r\to1-.$
\end{lem}

\begin{pf}
A simple computation gives
$$
\frac{zg'(z)}{g(z)}=1+\frac{\alpha z}{1-z}+
\frac {\beta z}{(1-z)\log\frac C{1-z}},
$$
where $C=e^{1/c}.$
By Lemma \ref{lem:C} and the fact that $\Re[z/(1-z)]>-1/2$ for $z\in\D,$
we obtain
$$
\Re\frac{zg'(z)}{g(z)}>
1-\frac\alpha2-\frac\beta{2\log(C/2)}
=1-\frac\alpha2-\frac\beta{\frac2c-2\log2}>0.
$$
Thus $g$ is starlike.
Since
$$\arg g_0(e^{i\theta})\to\pm\pi\alpha/2$$
as $\theta\to0\pm,$ we have $A(g)=\pi\alpha.$ The last assertion is
obvious.
\end{pf}

We are now in a position to prove Theorem \ref{thm:ce}.

\begin{pf}[Proof of Theorem \ref{thm:ce}]
Let $\lambda$ and $A$ be as in the theorem
and let $g$ be a function given in Lemma \ref{lem:example} for
$\alpha=A/\pi.$
We now define a function $f$ by the relation \eqref{eq:fg}.
Then $f\in\Sp_\lambda$ and, by \eqref{eq:M},
$$
M(r,f)\asymp M(r,g)^{\cos^2\lambda}
\asymp(1-r)^{-\alpha\cos^2\lambda}\big(\log\tfrac1{1-r}\big)^{\beta\cos^2\lambda}
$$
as $r\to1-.$
Therefore,
$$M(r,f)=O[(1-r)^{-A\cos^2\lambda/\pi}]$$
does not hold.
\end{pf}

\def\cprime{$'$} \def\cprime{$'$} \def\cprime{$'$}
\providecommand{\bysame}{\leavevmode\hbox to3em{\hrulefill}\thinspace}
\providecommand{\MR}{\relax\ifhmode\unskip\space\fi MR }
\providecommand{\MRhref}[2]{%
  \href{http://www.ams.org/mathscinet-getitem?mr=#1}{#2}
}
\providecommand{\href}[2]{#2}

\end{document}